\documentclass[12pt]{article}

\usepackage[utf8]{inputenc}

\usepackage[russian]{babel}

\usepackage{geometry}
\usepackage{color}
\usepackage{amsmath}
\usepackage{amssymb}
\usepackage{amsthm}

\newlength{\wdth}

\usepackage{hyperref}

\newtheorem{rDefinition}{Определение}
\newtheorem{rCorollary}{Следствие}
\newtheorem{rtheorem}{Теорема}
\newtheorem{rproposition}{Предложение}

\newif\ifeng
\engtrue 

\newif\ifhow
\howtrue 

\usepackage{hyperref}

\title{О Законе Больших Чисел 
для неодинаково распределенных слабо зависимых слагаемых
 \\~ \vspace{-2mm}\\
А.Т.Ахмярова\footnote{МГУ им. М.В.Ломоносова \& Институт проблем передачи информации РАН им. А.А. Харкевича, имейл: lili-r01@yandex.ru},   А.Ю.Веретенников\footnote{Институт проблем передачи информации РАН им. А.А. Харкевича \& Российский университет дружбы народов, имейл: ayv@iitp.ru} 

}

\begin{document}

\maketitle

\ifhow

\begin{abstract}
В работе предложены новые версии слабого Закона Больших Чисел для слабо зависимых слагаемых, вообще говоря, разнораспределенных, как при наличии математического ожидания каждого из них, так и без такового. Одним из основных условий для первого случая, который развивает идеи из статьи Y.S.Chow 1971г.,   являются равномерная интегрируемость слагаемых по Чезаро, в духе работ по ЗБЧ для попарно независимых случайных величин T.K.Chandra 1989 -- 2012г.г. При этом вместо попарной независимости налагаются совершенно иные условия слабой зависимости в духе статьи А.Н.Колмогорова 1929г., с той разницей, что в настоящей работе используются условия не на вторые, а лишь на первые моменты некоторых условных математических ожиданий. Второй случай основан на несколько ином условии слабой зависимости и использует телескопический метод и интерпретацию сходимости по вероятности к постоянной как слабую сходимость. Третий случай устанавливает ЗБЧ в отсутствие конечных математических ожиданий, опять же для разнораспределенных слагаемых.

\medskip

\noindent
Ключевые слова: Закон больших чисел; конечные математические ожидания; бесконечные математические ожидания;  слабая зависимость; разнораспределенные слагаемые.

\medskip


\end{abstract}

\section{Введение}
Пусть $\xi_1, \ldots, \xi_n, \ldots$ -- последовательность случайных величин на некотором вероятностном пространстве $(\Omega, {\cal F}, \mathsf P)$ и $S_n:= \sum_{i=1}^n \xi_i$. Законом Больших Чисел (далее ЗБЧ) называют утверждение вида  $S_n/n \stackrel{\mathsf P}\to 0$, или $(S_n - b)/n \stackrel{\mathsf P}\to 0$, или $(S_n - b_n)/n \stackrel{\mathsf P}\to 0$ с некоторыми  неслучайными $b$ или $b_n$. Третий вариант носит еще название устойчивости частот, см. \cite{Kolm}. Также все эти варианты  часто называют слабым ЗБЧ, чтобы отличить от усиленного ЗБЧ, в котором сходимость имеет место почти наверное (см. \cite{Kolm}).

Первый строго доказанный вариант ЗБЧ появился в опубликованном в 1713 году труде Я.Бернулли (Jakob Bernoulli) \cite{Bern} для схемы Бернулли (разумеется, тут цитируется не столь давний перевод на русский, а не первоисточник). Ранее  в этой области  высказывались лишь нестрогие утверждения. Важнейшим открытием 19 века явилась работа П.Л. Чебышева \cite{Chebyshev}. О дальнейшем, как и о предшествующем  прогрессе в этой тематике написать кратко не представляется возможным, поэтому отошлем читателя к подробным обзорам \cite{Seneta13} и \cite{Ash21}, а также к монографии \cite{Chandra}. 
Мотивация не требует подробных объяснений, поскольку закон  больших чисел является признанной базовой теоремой теории вероятностей и математической статистики. 



В этой работе предложены новые варианты (слабого) ЗБЧ для слабо зависимых слагаемых, вообще говоря, разнораспределенных, как с конечными, так и с бесконечными математическими ожиданиями.

Скажем несколько слов об условиях слабой зависимости, при которых известен ЗБЧ. Исторически первое такое условие --  полная независимость, то есть, независимость всех слагаемых в совокупности. Многие современные варианты ЗБЧ предполагают лишь попарную независимость, в том числе, в \cite{Etemadi}. 
Версии ЗБЧ также известны для эргодических марковских процессов с не слишком медленным перемешиванием, а для конечных цепей Маркова с экспоненциальной скоростью сходимости эта теорема была  доказана еще самим А.А.Марковым \cite{Markov}; в принципе, хорошее перемешивание может вести к ЗБЧ и без условия связанности в цепь Маркова, однако,  перемешивание не является темой данной работы. Здесь также уместно упомянуть об эргодической теореме Биркгофа -- Хинчина для лишь стационарных случайных величин с конечным математическим ожиданием, см. \cite{Kolm_erg}, хотя в строгом смысле этот результат не является законом больших чисел, поскольку в нем сходимость устанавливается, вообще говоря, к некоторой случайной величине с, вообще говоря, неизвестным распределением.

\medskip

Далее в настоящей работе будут использованы условия слабой зависимости такого типа: при условии $\mathsf E\xi_n=0$ предполагается, что $\mathsf E(\xi_n|\xi_1, \ldots, \xi_{n-1}) \underset{n\to\infty}{\to} 0$, либо $\mathsf E(|\xi_n|\, |\xi_1+\ldots + \xi_{n-1}) \underset{n\to\infty}{\to} 0$, в том числе,  в форме ``по Чезаро'' (точные формулировки см. в следующих разделах). А.Н.Колмогоров указывал похожие (только без Чезаро), хотя формально несколько иные условия \cite[гл.6.3]{Kolm}, \cite{Kolm1929}, связанные со вторыми моментами некоторых условных математических ожиданий. В случае бесконечных математических ожиданий будет использован аналог вида
\(
n^{-1}\sum_{k=1}^n \mathsf E(\xi_k 1(|\xi_k|\le n) | \xi_1, \ldots, \xi_{k-1}) \underset{n\to\infty}{\to} 0
\).

\medskip




Работа состоит из трех разделов: данное введение, основные результаты -- теоремы 1--3 -- и доказательства, разбитые на подразделы. Ради удобства читателя вместе с основными результатами приведены некоторые известные классические варианты теорем.   Все новые результаты обозначаются как теоремы, а все классические -- как предложения. Теорема \ref{rthm1} является основным результатом курсовой работы первого автора. Теоремы \ref{rthm2} и \ref{rthm3} установлены вторым автором.

\section{Основные результаты}
Прежде всего, напомним определения равномерной интегрируемости (далее Р.И.) и Р.И. по Чезаро. 

\begin{rDefinition}
Последовательность с.в. $(\xi_n)$ называется Р.И., если 
$$
\lim_{M\to\infty} \sup_n \mathsf E |\xi_n| 1(|\xi_n| >M) = 0. 
$$
\end{rDefinition}

\begin{rproposition}[де Ла Валле Пуссен]
Последовательность Р.И. тогда и только тогда, когда найдется функция $g(t)$ такая, что $g(t)/t \uparrow \infty, t\to\infty$, $g$ выпукла и  $$\sup_n \mathsf E g(|\xi_n|)<\infty.$$
\end{rproposition}

\begin{rDefinition}
Последовательность с.в. $(\xi_n)$ называется Р.И. по Чезаро, если 
$$
\lim_{M\to\infty} \sup_n \frac1{n} \sum_{k=1}^n \mathsf E |\xi_k| 1(|\xi_k| >M) = 0. 
$$
\end{rDefinition}

\subsection{Теорема \ref{rthm1}}

Первая теорема является развитием, или обобщением классических результатов Хинчина (1929) и Y.S.Chow (1971). Приведем сперва их формулировки.

\begin{rproposition}[А.Я.Хинчин \cite{Khintchine}]
Пусть $\xi_1, \ldots, \xi_n$ -- независимые, одинаково распределенные случайные величины с конечным математическим ожиданием $\mathsf E\xi_1 = 0$. Тогда 
$S_n/n \stackrel{\mathsf P}\to 0$. 
\end{rproposition}
Стоит напомнить, что в тех же предположениях имеет место теорема Колмогорова -- усиленный ЗБЧ, см. \cite[Теорема VI.5.III]{Kolm}.

Хотя результат N.Etemadi относится к усиленному ЗБЧ (УЗБЧ), здесь также будет уместно его напомнить, поскольку, как хорошо известно, УЗБЧ влечет справедливость и слабого ЗБЧ. Полезно еще отметить, как это подчеркнуто в \cite{Chibisov}, что слабый ЗБЧ допускает несколько более общую формулировку, чем в большинстве работ из этой области, а именно, в виде схемы серий, когда каждая серия может быть определена на своем вероятностном пространстве. В данной работе этот вариант не обсуждается, однако, его надо иметь в виду.

\begin{rproposition}[N.Etemadi \cite{Etemadi}]
Пусть $(\xi_n, n\ge 1)$ -- попарно независимые, одинаково распределенные с.в. с конечным математическим ожиданием, $\mathsf E \xi_n = 0, \, \forall \, n$. Тогда  $S_n/n\stackrel{\text{\small п.н. \& }\mathsf P}\to 0, \; n\to\infty$.

\end{rproposition}

Образцом для нашей первой теоремы послужит результат Чау, который будет здесь приведен только для случая сходимости в $L_p$ при $p=1$; в \cite{Chow} рассмотрено $0<p<2$.
\begin{rproposition}[Y.S.Chow \cite{Chow}]\label{Chow}
Пусть семейство случайных величин $(\xi_n, n\ge 1)$ Р.И., все они имеют конечное нулевое математическое ожидание $\mathsf E \xi_n = 0$, и 
$$
a_n:= \sum_{k=1}^n \mathsf E(\xi_k|\xi_1, \ldots, \xi_{k-1}).
$$ 
Тогда
$$
\mathsf E |S_n/n - a_n/n|\to 0, \quad n\to\infty.
$$
\end{rproposition}
Строго говоря, это еще не ЗБЧ, хотя и очень близко. По какой-то причине автор не сформулировал данное очевидное следствие из своего результата. Сделаем это за него (конечно, не претендуя на авторство). 

\begin{rCorollary}[ЗБЧ при слабой зависимости $a_n/n \to 0$]
Если в условиях предложения \ref{Chow} выполнено условие слабой зависимости
$$
a_n/n \stackrel{L_1}\to 0, \quad n\to\infty,
$$ 
то 
$$
\mathsf E |S_n/n|\to 0, \quad n\to\infty.
$$
Если же вместо этого лишь
$
a_n/n \stackrel{\mathsf P}\to 0$, $n\to\infty, 
$
то 
$$
S_n/n \stackrel{\mathsf P}\to 0, \quad n\to\infty.
$$
\end{rCorollary}

Классический результат Чандры относится к попарно независимым слагаемым. 
\begin{rproposition}[T.K.Chandra  \cite{Chandra}]
Пусть $(\xi_n, n\ge 1)$ -- попарно независимые не о.р. с.в. Р.И. по Чезаро, $\mathsf E \xi_n = 0, \, \forall \, n$. Тогда  $S_n/n\stackrel{\mathsf P \& L_1}\to 0, \; n\to\infty$.
\end{rproposition}

Следующая теорема является основным результатом курсовой работы первого автора 2023--2024г. и одновременно первым основным результатом данной работы.

\begin{rtheorem}[]\label{rthm1}
1. Пусть семейство случайных величин $(\xi_n, n\ge 1)$ Р.И. по Чезаро, все они имеют конечное нулевое математическое ожидание $\mathsf E \xi_n = 0$, и $a_n:= \sum_{k=1}^n \mathsf E(\xi_k|\xi_1, \ldots, \xi_{k-1})$. Тогда
\begin{equation}\label{chow1}
\mathsf E |S_n/n - a_n/n|\to 0, \quad n\to\infty.
\end{equation}
2. Если при тех же условиях еще 
\begin{equation}\label{wd1}
a_n/n \stackrel{\mathsf P}\to 0, \quad n\to\infty,
\end{equation} 
то 
\begin{equation}\label{lln1}
S_n/n \stackrel{\mathsf P}\to 0, \quad n\to\infty.
\end{equation} 
\end{rtheorem}
\noindent
Здесь условие (\ref{wd1}) можно считать условием слабой зависимости.


\subsection{Теорема \ref{rthm2}}

В следующем результате используется несколько иное свойство слабой зависимости: условные математические ожидания для каждой случайной величины вычисляются при условии {\em суммы} всех предыдущих слагаемых.
\begin{rtheorem}\label{rthm2}
Пусть семейство случайных величин $(\xi_n, n\ge 1)$ Р.И. по Чезаро, все они имеют конечные нулевые математические ожидания $\mathsf E \xi_n = 0$, $n\ge 1$, и выполнено условие слабой зависимости
\begin{equation}\label{CUI2}
 \frac1 {n} \sum_{k=1}^n \mathsf E |\mathsf E(\xi_{k} |\xi_{1}+ \ldots+ \xi_{k-1})|  \to 0, \quad n\to\infty.
\end{equation}
Тогда
$$
S_n/n \stackrel{\mathsf P}\to 0, \quad n\to\infty.
$$
\end{rtheorem}
Отметим, что достаточным условием для справедливости  (\ref{CUI2})  (по модулю  Р.И. по Чезаро) является сходимость 
$$
 \mathsf E(\xi_{k} |\xi_{1}+ \ldots+ \xi_{k-1})\stackrel{\text{a.s.}} \to 0, \quad k\to\infty, 
$$
или
$$
 \mathsf E(\xi_{k} |\xi_{1}+ \ldots+ \xi_{k-1})\stackrel{\mathsf P} \to 0, \quad k\to\infty.
$$
Отметим также, что условие Р.И. по Чезаро влечет за собой оценку
$$
\sup_n \frac1{n} \sum_{k=1}^n \mathsf E |\xi_k| <\infty.
$$

\medskip


\subsection{Теорема \ref{rthm3}}
Для формулировки следующих результатов напомним обозначение $\bar F(x):=1-F(x)$, где $F(x)$ -- произвольная функция распределения. Подчеркнем, что далее существование математических ожиданий $\xi_n$ не предполагается, хотя возможно их  существование в обобщенном смысле ``А-интеграла'', то есть, в смысле главного значения (введено Титчмаршем).

\begin{rproposition}[Колмогоров, \cite{Kolm1929}, также гл. VI, \S 4 \cite{Kolm}]\label{pro6}
Пусть $(\xi_n, n\ge 1)$ -- независимые, одинаково распределенные случайные величины и выполнено условие 
\begin{equation}\label{ekolm}
\lim\limits_{n\to \infty}n(\bar F_\xi(n) + F_\xi(-n)) =  0. 
\end{equation}
Тогда (и только тогда) $S_n/n - \mu_n \stackrel{\mathsf P}\to 0, \; n\to\infty$,  где $\mu_n= \mathsf E \xi_1 1(|\xi_1|\le n)$.
\end{rproposition}

Замечание. Возможна ситуация, когда существует нулевое математическое ожидание в смысле А-интеграла (главного значения), то есть, $\mu_n\to 0, \, n\to\infty$. Например, для симметричных распределений просто $\mu_n = 0$ при всех $n$. В этом случае (то есть, при $\mu_n\to 0, \, n\to\infty$) в условиях предложения \ref{pro6} имеем $S_n/n \stackrel{\mathsf P}\to 0, \; n\to\infty$.



Здесь,  $n(\bar F_{\xi_1}(n) + F_{\xi_1}(-n))$ оценивает сверху вероятность $\tau_n=\mathsf P(\bigcup_{k=1}^n (|\xi_k|>n))$, так что условие (\ref{ekolm}) в случае независимых, одинаково распределенных слагаемых  позволяет применить метод урезания на уровне $n$. Обозначим
$$
\gamma_k(x)=\bar F_{\xi_k}(x) + F_{\xi_k}(-x).
$$ 
В общем случае когда равенство распределений $\xi_n$ не предполагается, будем оценивать вероятность $\tau_n$ как $\tau_n\le T_n:= \sum_{k=1}^n (\bar F_{\xi_k}( n) + F_{\xi_k}(-n))=\sum_{k=1}^n\gamma_k(n)$. 
Пусть 
\[
\tilde a_n/n := n^{-1}\sum_{k=1}^n \mathsf E(\xi_k 1(|\xi_k|\le n) | \xi_1, \ldots, \xi_{k-1}).
\]
Положим также 
$$
\psi_n(y):=\sum\limits_{k=1}^{n} y\gamma_k(ny).
$$
Следующая теорема является третьим основным результатом данной работы. 
\begin{rtheorem}[для разнораспределенных с.в. без математических ожиданий]\label{rthm3}
Пусть случайные величины $(\xi_n, n\ge 1)$  
таковы, что 
семейство функций $(\psi_n(y), \, 0 \le y \le 1)$ является равномерно интегрируемым, а также при всех $y\in [0,1]$
\begin{equation}\label{psi0}
\lim\limits_{n\to \infty}\psi_n(y) =  0,
\end{equation}
и выполнено условие слабой зависимости 
\begin{equation}\label{wd30}
\tilde a_n/n \overset{\mathsf P}{\underset{n\to\infty}{\to}} 0.
\end{equation}
Тогда 
\begin{equation}\label{lln3}
S_n/n  \stackrel{\mathsf P} {\underset{n\to\infty}{\to}} 0.
\end{equation}  
\end{rtheorem}
Напомним, что {\em достаточными} условиями для равномерной интегрируемости семейства функций $\left(\sum\limits_{k=1}^{n} y\gamma_k(ny), \, 0 \le y \le 1\right)$ являются их равностепенная ограниченность, либо интегрируемость вида $\displaystyle \int_0^1 \sup_n \sum\limits_{k=1}^{n} y\gamma_k(ny) dy <\infty$. Также отметим, что условие (\ref{psi0}) в случае независимых одинаково распределенных слагаемых превращается в (\ref{ekolm}), а функция в левой части (\ref{ekolm}) автоматически является ограниченной при выполнении этого условия. Таким образом, можно считать, что условие (\ref{psi0}) является адекватным аналогом (\ref{ekolm}) для случая неодинаково распределенных и даже не обязательно независимых случайных величин. Условие (\ref{wd30}) можно считать вариантом слабой зависимости $(\xi_k)$ по аналогии с условиями (\ref{wd1}) и (\ref{CUI2}), хотя в данном случае конечность математических ожиданий $\xi_k$ не предполагается. В случае независимости в совокупности всех $\xi_k$ условие (\ref{wd30}) превращается в следующее соотношение: 
\begin{align}\label{nrv}
\frac1n \sum_{k=1}^n \mathsf E\xi_k 1(|\xi_k|\le n) \to 0, \quad n\to\infty, 
\end{align}
что является версией сходимости к нулю ``по Чезаро'' математических ожиданий первых $n$ слагаемых в смысле А-интеграла, то есть, главного значения. При этом, из-за добавления  ``по Чезаро'' ни одна из величин $\xi_k$ в отдельности ни при каком фиксированном $k$ не обязана иметь нулевое, или близкое к нулевому математическое ожидание в смысле главного значения. Даже в таком варианте с независимыми в совокупности, но не одинаково распределенными $\xi_k$ данный результат, по-видимому, является новым.

\section{Доказательства}
\subsection{Доказательство теоремы 1}
Пусть $M>0$. Положим 
\[
\xi'_n=\xi_n1[|\xi_n|\leq M], \quad \xi''_n=\xi_n-\xi'_n.
\]
{\bf 1}. Покажем, что ряд \(\sum_{k=1}^\infty {k^{-1} (\xi'_k- \mathsf E(\xi'_k|\xi_1,...,\xi_{k-1}))}\) сходится.
При $j<i$ имеем:
\begin{align*}
&\mathsf E[ (\xi'_i-E(\xi'_i|{\cal F}^\xi_{i-1}))(\xi'_j-\mathsf E(\xi'_j|{\cal F}^\xi_{j-1}))]
\\\\
&=\mathsf  E (\xi'_j- \mathsf E(\xi'_j|{\cal F}^\xi_{j-1})\mathsf E[ (\xi'_i- \mathsf E(\xi'_i|{\cal F}^\xi_{i-1}))|{\cal F}^\xi_{i-1}]
\\\\
&= \mathsf E (\xi'_j-E(\xi'_j|{\cal F}^\xi_{j-1}) [\mathsf E (\xi'_i|{\cal F}^\xi_{i-1})) - \mathsf E (\xi'_i|{\cal F}^\xi_{i-1}))] = 0. 
\end{align*}
При $j=i$ имеем оценку:
\[
 \mathsf E(\xi'_i-  \mathsf  E(\xi'_i|{\cal F}^\xi_{i-1}))^2=  \mathsf E[2((\xi'_i)^2+2( \mathsf E(\xi'_i|{\cal F}^\xi_{i-1}))^2] \leq 2M^2+2M^2 \leq 4M^2.
\]
Складывая, получаем 
\[ 
\sum_{k=1}^\infty  \mathsf  D(\xi'_k-E(\xi'_k|{\cal F}^\xi_{k-1}))/k^2 \leq 4M^2/k^2 < \infty. 
\]
Также заметим, что
\begin{align*}
& \mathsf E \sum_{k=1}^n {k^{-1} (\xi'_k- \mathsf E(\xi'_k|\xi_1,...\xi_{k-1}))}
=\sum_{k=1}^n {\frac1k  \mathsf E(\xi'_k- \mathsf E(\xi'_k|\xi_1,...\xi_{k-1}))}
=0.
\end{align*}
Согласно 
теореме Хинчина - Колмогорова о сходимости ряда из случайных величин \cite{Khin-Kolm}, доказательство которой остается верным при лишь некореллированности слагаемых, заключаем, что, как и было обещано, сходится ряд 
\[
\sum_{k=1}^\infty {k^{-1} (\xi'_k- \mathsf E(\xi'_k|\xi_1,...\xi_{k-1}))}<\infty\quad \text{п.н.}
\]

\medskip

\noindent
{\bf 2}. Теперь в силу леммы Кронекера имеем, 
\[ 
\frac1n \sum_{k=1}^n [\xi'_k- \mathsf E(\xi'_k|\xi_1,...\xi_{k-1}))] =o(1), \quad n\to\infty.
\]

\noindent
{\bf 3}. В силу предположения Р.И. по Чезаро, выполнены условия  теоремы Витали о предельном переходе под знаком интеграла. Стало быть, заключаем, что
\[ 
\lim_{n\to\infty} \mathsf E\left|\frac1 n \sum_{k=1}^n [\xi'_k- \mathsf E(\xi'_k|\xi_1,...\xi_{k-1}))]\right| =0.
\]

\noindent
{\bf 4}. Далее, оцениваем,
\begin{align*}
& \mathsf E|\sum_1^n [\xi''_k- \mathsf E(\xi''_k|\xi_1,...\xi_{k-1})]| \leq 2 \mathsf E \sum_1^n|\xi''_k|
=2 \mathsf E \sum_1^n|\xi_k-\xi'_k| 
\\
&= 2 \mathsf E \sum_1^n|\xi_k 1[|\xi_k|>M]|
=2 \mathsf E \sum_1^n|\xi_k| 1[|\xi_k|>M] \leq 2n\epsilon
\end{align*}

Вновь в силу условия Р.И. по Чезаро, для всякого $\epsilon >0$ найдется такое $M_0$, что при всех $M\geq M_0$, 
\[
\frac 1 n  \sum_1^n  \mathsf E|\xi_k|1[|\xi_k|>M]<\epsilon.
\]

\noindent
{\bf 5}. В итоге получаем, 
\begin{align*}
&\frac1n \mathsf E|\sum_1^n \xi_k-\sum_1^n \mathsf E(\xi_k|\xi_1,...\xi_{k-1})|
\\\\
&= \mathsf E\frac 1 n| \sum_1^n \left[ \xi'_k- \mathsf E(\xi'_k|\xi_1,...\xi_{k-1})+\xi''_k- \mathsf E(\xi''_k|\xi_1,...\xi_{k-1}) \right]| 
\\\\
&\leq o(1)+2 \epsilon, \quad n\to\infty.
\end{align*}
Здесь левая часть не зависит от $\epsilon$. Поэтому, 
устремляя $\epsilon$ к 0, получаем
\[
\frac1n \mathsf E|S_n-a_n|\underset{n \to \infty}{\longrightarrow} 0,
\]
где $S_n=\xi_1+...+\xi_n$, $a_n=\sum_1^n \mathsf E(\xi_k|\xi_1,...\xi_{k-1})$, что и требовалось показать. Первое утверждение теоремы 1 доказана. 

\medskip

\noindent
Второе утверждение -- собственно ЗБЧ (\ref{lln1}) -- при условии слабой зависимости (\ref{wd1}) следует из (\ref{chow1}) непосредственно. \hfill QED

\fi

\subsection{Доказательство теоремы 2}
Используем метод ``телескопического разложения''. Как хорошо известно, он с успехом применяется в доказательстве ЦПТ. Для доказательства же ЗБЧ такой метод был, возможно, впервые  предложен в учебном пособии автора \cite{Ver_MIREA-2} (хотя и с некоторыми опечатками) при условии независимости и одинаковой распределенности слагаемых. Как оказалось, он работает и в условиях слабой зависимости (\ref{CUI2}).

\noindent 
{\bf 0.} Как хорошо известно, сходимость по вероятности {\bf к константе} $S_n/n \stackrel{\mathsf P}\to 0$ может быть интерпретирована как слабая сходимость:
\begin{equation}\label{weak}
\mathsf E f(S_n/n) \to f(0), \quad n\to\infty, 
\end{equation}
для любой функции $f\in C_b(\mathbb R)$. Будем доказывать указанную слабую сходимость. 

\medskip

\noindent
{\bf 1.} Еще отметим, что для того, чтобы установить соотношение  (\ref{weak}), достаточно его проверить для более узкого класса функций, в частности, для любой $f\in C^1_b(\mathbb R)$ с ограниченным модулем непрерывности производной $\rho = \rho_{f'}$. В самом деле, семейство таких функций всюду плотно в пространстве $C_b([-N,N])$ для любого $N$, что позволяет аппроксимировать и затем переходить к пределу. (В частности,  достаточно установить такую сходимость лишь для функций вида $f(x) = \sin(c x +\varphi_0)$; на этом замечании основан подход к предельным теоремам на основе характеристических функций.)

\noindent
{\bf 2.} Для любой измеримой ограниченной функции $f$ имеем с  $S_0=0$, 
\begin{align}\label{fn0}
\mathsf E f(S_n/n) - f(0) = \sum_{k=1}^{n} (\mathsf E f(S_k/n) - \mathsf E f(S_{k-1}/n) ).
\end{align}
Согласно версии теоремы Ньютона -- Лейбница при $f\in C^1_b$, 
$$
f(y) - f(x) = (y-x)\int_0^1 f'(x+a(y-x)) da, 
$$
имеем, 
\begin{align*}
f(S_k/n) - f(S_{k-1}/n) 
= \frac{\xi_k}{n} \int_0^1 f'(S_{k-1}/n + a \xi_k/n) da.
\end{align*}

\noindent
{\bf 3.} 
Под  интегралом хотелось бы вычесть $f'(S_{k-1}/n)$, чтобы получить $\int_0^1 (f'(S_{k-1}/n + a \xi_k/n) -f'(S_{k-1}/n)) da$.
С этой целью заметим, что в силу условия $\|f'\| < \infty$ справедлива оценка
\begin{align*}
&|\mathsf E \frac{\xi_k}{n} \int_0^1 f'(S_{k-1}/n) da |
= \frac1{n} |\mathsf E \int_0^1 f'(S_{k-1}/n) da \, \mathsf E (\xi_k | S_{k-1})| 
 \\&
\le \frac{\|f'\| }{n}  \mathsf E |\mathsf E (\xi_k | S_{k-1})|.
\end{align*}
Поэтому, в силу предположения теоремы (\ref{CUI2}), 
\begin{align}\label{r12}
&\sum_{k=1}^n |\mathsf E \frac{\xi_k}{n} \int_0^1 f'(S_{k-1}/n) da |
= \frac1{n} \sum_{k=1}^n |\mathsf E \int_0^1 f'(S_{k-1}/n) da \mathsf E (\xi_k | S_{k-1})| 
\nonumber
 \\ 
&\le  \frac{\|f'\|}{n}  \sum_{k=1}^n  \mathsf E |\mathsf E (\xi_k | S_{k-1})| \to 0, \quad n\to\infty.
\end{align}

\noindent
{\bf 4.}
Итак, оцениваем 
\begin{align}
&| \sum_{k=1}^n (\mathsf E f(S_k/n) - \mathsf E   f(S_{k-1}/n))| 
= |\sum_{k=1}^n \mathsf E \frac{\xi_k}{n} \int_0^1 f'(S_{k-1}/n + a \xi_k/n) da|
\nonumber
 \\ \nonumber\\
&=  |\sum_{k=1}^n \mathsf E \frac{\xi_k}{n} \int_0^1 (f'(S_{k-1}/n + a \xi_k/n) - f'(S_{k-1}/n)) da| + o(1)
\nonumber
 \\ \nonumber\\
\label{o1}
&\le
 \sum_{k=1}^n \mathsf E \frac{|\xi_k|}{n} \rho_{f'}(\xi_k/n)  + o(1)
 = \mathsf E \sum_{k=1}^n \frac{|\xi_k|}{n} \rho_{f'}(\xi_k/n)  + o(1).
\end{align}

\noindent
{\bf 5.} Покажем, что последнее математическое ожидание здесь также стремится к нулю при $n \to \infty$.  Напомним, что функция $\rho_{f'}$ ограничена. Оценим величину 
$\mathsf E \sum_{k=1}^n \frac{|\xi_k|}{n} \rho_{f'}(\xi_k/n)$\; при условии $\sup_n \frac1{n} \sum_{k=1}^n \mathsf E |\xi_k| 1(|\xi_k| >M) {\underset{M\to\infty}{\to}} 0$. 
Пусть $\varepsilon >0$. 
Имеем, при всяком $M>0$, 
\begin{align*}
& \frac{1}{n} \sum_{k=1}^n  \mathsf E |\xi_k| \rho_{f'}(\xi_k/n)
=  \frac{1}{n} \sum_{k=1}^n \mathsf E |\xi_k|\rho_{f'}(\xi_k/n) 1(|\xi_k|\le M)
 \\\\
& +\frac{1}{n} \sum_{k=1}^n \mathsf E |\xi_k| 1(|\xi_k| > M) \underbrace{\rho_{f'}(\xi_k/n)}_{\le K}=: \Sigma^1_n(M) + \Sigma ^2_n(M).
\end{align*}
Согласно условию Чезаро Р.И., можно выбрать такое $M$, что $\Sigma ^2_n(M) \le \varepsilon$. При таком $M$ имеем оценку для $\Sigma ^1_n(M)$:
$$
\Sigma ^1_n(M) \le M  \rho_{f'}(M/n) \to 0, \quad n\to \infty,
$$
в силу непрерывности $\rho_{f'}$ в нуле и равенства $\rho_{f'}(0)=0$.
Это и доказывает, что выражение в правой части (\ref{o1}) стремится к нулю при $ n\to \infty$, что и требовалось. \hfill QED

\subsection{Доказательство теоремы \ref{rthm3}}
Воспользуемся урезанием на уровне $n$, то есть, положим 
$$
\xi'_{k,n}=\xi'_k := \xi_k 1(|\xi_k|\le n), \quad 
S'_n:=\sum_{k=1}^n \xi'_{k}.
$$
Согласно определению случайных величин $\tilde a_n$, 
\[
\tilde a_n/n 
= n^{-1}\sum_{k=1}^n \mathsf E(\xi'_k| \xi_1, \ldots, \xi_{k-1}).
\]
Оценим сначала вероятность $\mathsf P(\frac1n |S_n - \tilde a_n|>x)$. 
Как легко видеть, и этот прием применялся в литературе, начиная, видимо, с работы \cite{Khintchine}, 
\begin{equation}\label{trunc}
\mathsf P(\frac1n |S_n - \tilde a_n|>x) \le \mathsf P(|S'_n - \tilde a_n|>nx) + \mathsf P(S_n \neq S'_n).
\end{equation}
Имеем, 
\begin{equation}\label{sneqs}
\mathsf P(S_n \neq S'_n) \le \sum_{k=1}^n \mathsf P(|\xi_{k}|>n) 
= \sum_{k=1}^n \gamma_k(n) = T_n = \psi_n(1) \to 0, \quad n\to\infty, 
\end{equation}
по условию теоремы. 

\medskip

Оценим теперь первое слагаемое в правой части (\ref{trunc}). Запишем,
\begin{align*}
S'_n - \tilde a_n = 
\sum_{k=1}^n (\xi'_{k} - \mathsf E(\xi'_k | \xi_1, \ldots, \xi_{k-1})).
\end{align*}
Здесь
\begin{align*}
\mathsf E(S'_n - \tilde a_n) = 0, 
\end{align*} 
поэтому далее оценим дисперсию случайной величины $(S'_n - \tilde a_n)/n$. В силу некоррелированности случайных величин $\xi'_{k} - \mathsf E(\xi'_k | \xi_1, \ldots, \xi_{k-1})$ при различных $k$, которая доказывается точно так же\footnote{Имеем при $i<j$ (эту выкладку авторы предлагают не включать в текст, если работа будет принята: здесь она приведена только ради небольшого облегчения чтения рецензента), 
\begin{align*}
&\mathsf E (\xi'_i - \mathsf E(\xi'_i|{\cal F}_{i-1}) (\xi'_j - \mathsf E(\xi'_j|{\cal F}_{j-1}) 
= \mathsf E \mathsf E[(\xi'_i - \mathsf E(\xi'_i|{\cal F}_{i-1}) (\xi'_j - \mathsf E(\xi'_j|{\cal F}_{j-1})|{\cal F}_{j-1}]
 \\\\
&= \mathsf E (\xi'_i -\mathsf  E(\xi'_i|{\cal F}_{i-1})  \mathsf E[(\xi'_j - \mathsf E(\xi'_j|{\cal F}_{j-1})|{\cal F}_{j-1}] 
=\mathsf  E (\xi'_i - \mathsf E(\xi'_i|{\cal F}_{i-1}) \{E(\xi'_j|{\cal F}_{j-1}) - \mathsf E(\xi'_j|{\cal F}_{j-1})\} =0. 
\end{align*} }, как в п. 1 доказательства теоремы \ref{rthm1}, находим,
\begin{align*}
&\mathsf D((S'_n - \tilde a_n)/n) = 
\frac1{n^2} \!\sum_{k=1}^n \mathsf D(\xi'_k \!-\! \mathsf E(\xi'_k|{\cal F}^\xi_{k-1}))  
=  \frac1{n^2} \!\sum_{k=1}^n \!\mathsf E (\xi'_k - \mathsf E(\xi'_k|{\cal F}^\xi_{k-1}))^2  
 \\\\
&\le  \frac2{n^2} \sum_{k=1}^n (\mathsf E (\xi'_k)^2 + \mathsf E(\mathsf E(\xi'_k|{\cal F}^\xi_{k-1}))^2) 
\le  \frac2{n^2} \sum_{k=1}^n (\mathsf E (\xi'_k)^2 + \mathsf E(\mathsf E((\xi'_k)^2|{\cal F}^\xi_{k-1})))
 \\\\
&=  \frac2{n^2} \sum_{k=1}^n (\mathsf E (\xi'_k)^2 
+ \mathsf E(\xi'_k)^2)) 
= \frac4{n^2} \sum_{k=1}^n \mathsf E (\xi'_k)^2.
\end{align*}
Следуя \cite[Т. 2, гл. VII, \S 7, формула (7.7)]{Feller},
введем обозначение 
$$
\sigma_k(t) : = \frac1t \int_{-t}^t x^2 dF_{\xi_k}(x) = - t \gamma_k(t) + \frac2t \int_0^t x\gamma_k(x)dx.
$$
Отметим, что в \cite{Feller} разобран лишь случай независимых, одинаково распределенных величин; в нашем же случае мы обязаны добавлять индекс $k$ у всех слагаемых $\xi_k$; стало быть, и величина $\sigma_k(t)$ должна иметь этот же индекс. Второе равенство здесь следует из интегрирования по частям; данная формула исправлена по сравнению с ошибочной формулировкой 
\cite[Т. 2, гл. VII, (7.7)]{Feller}; отметим, однако, что правильное {\em неравенство,} вытекающее из этой формулы и используемое далее для доказательство ЗБЧ, можно найти в \cite[Т. 2, гл.  XVII, формула (2.39)]{Feller}. В частности, 
$$
\frac1n \mathsf E (\xi'_k)^2 = \sigma_k(n) : = \frac1n \int_{-n}^n x^2 dF_{\xi_k}(x) = - n \gamma_k(n) + \frac2n \int_0^n x\gamma_k(x)dx.
$$
Поскольку $\gamma_k(n)\ge 0$, то, отбрасывая первый член в последнем равенстве, в силу условий теоремы о равномерной интегрируемости и сходимости к нулю на $\psi_n$ получаем, 
\begin{align*}
\frac4{n^2} \sum_{k=1}^n \mathsf E (\xi'_k)^2
\le  \frac8{n^2} \int\limits_0^n \sum_{k=1}^n x\gamma_k(x)dx  \stackrel{y=x/n}= 8 \int\limits_0^1 \underbrace{\sum_{k=1}^n y\gamma_k(ny)}_{=\psi_n(y)}dy \to 0, \quad n\to\infty.
\end{align*}
Стало быть, в силу неравенства Бьенаме -- Чебышева,
\begin{align}\label{sprima}
&\mathsf P(|S'_n - \tilde a_n|/n >\epsilon) \le \epsilon^{-2}\mathsf D((S'_n - \tilde a_n)/n) \to 0, \quad n\to\infty.
\end{align}

Из (\ref{trunc}), (\ref{sneqs}), (\ref{sprima}) и условия (\ref{wd30}), наконец, заключаем, что имеет место сходимость (\ref{lln3}), что и требовалось. \hfill QED

\section*{Благодарности}
Для обоих авторов работа поддержана Фондом развития теоретической физики и математики ``Базис''.


\end{document}